\documentclass[12pt,reqno]{amsart}
\usepackage{amsmath,amssymb,amsfonts}
\usepackage[dvips]{graphicx,color}
\usepackage{mathrsfs}
\usepackage[abbrev]{amsrefs}

\usepackage{latexsym}
\usepackage{delarray}

\numberwithin{equation}{section}
\theoremstyle{plain}
\newtheorem{thm}{Theorem}[section]
\newtheorem{lem}[thm]{Lemma}

\newtheorem{cor}[thm]{Corollary}
\theoremstyle{definition}
\newtheorem{defi}[thm]{Definition}
\newtheorem{rem}[thm]{Remark}

\newcommand\R{\mathbb{R}}

\numberwithin{equation}{section}
\numberwithin{thm}{section}

\makeatletter
\newcommand{\vertiii}[1]{{\left\vert\kern-0.25ex\left\vert\kern-0.25ex\left\vert #1
    \right\vert\kern-0.25ex\right\vert\kern-0.25ex\right\vert}}

\title[Gradient estimates   for  heat equations]{
Gradient estimates and their optimality for  heat equation in an exterior domain}
\author[V. Georgiev and
K. Taniguchi]{Vladimir Georgiev and Koichi Taniguchi}
\date{\today}

\address{Vladimir Georgiev \\ Department of Mathematics, University of Pisa, Largo B. Pontecorvo 5 Pisa,
56127 Italy \\ and \\ Faculty of Science and Engineering, Waseda University,
 3-4-1, Okubo, Shinjuku-ku, Tokyo 169-8555,
Japan
}
\email{georgiev@dm.unipi.it}

\address{Koichi Taniguchi \\ Department of Mathematics,
Chuo University,
1-13-27, Kasuga, Bunkyo-ku,
Tokyo 112-8551,
Japan
}
\email{koichi-t@gug.math.chuo-u.ac.jp}

\thanks{The first  author  was supported in part by  INDAM, GNAMPA - Gruppo Nazionale per l'Analisi Matematica, la Probabilita e le loro Applicazion, by Institute of Mathematics and Informatics, Bulgarian Academy of Sciences and by Top Global University Project, Waseda University}

\keywords{Heat equations, gradient estimates, Dirichlet problem, exterior domains}

\begin{document}


\footnote[0]
{2010 {\it Mathematics Subject Classification.}
Primary 35K05; Secondary 35K20.}

\begin{abstract}
This paper is devoted to the study of gradient estimates
for the Dirichlet problem of the
heat equation in the exterior domain of a compact set.
Our results describe the time decay rates of the derivatives of solutions to the Dirichlet problem.  Applications of these estimates to bilinear type commutator estimates for Laplace operator with Dirichlet boundary condition in exterior domain are discussed too.
\end{abstract}

\maketitle


\section{Introduction}

Studying gradient estimates over the evolution flow of parabolic equations is a challenging problem having different applications in the theory of incompressible Navier - Stokes flow (see \cite{DSh99}, \cite{Han-2012}, \cite{His-2016}, \cite{HisShi-2009}, \cite{IftKarLac-2014}, \cite{KobKub-2013}, \cite{MarSol-1997})
as well in  harmonic analysis (comparison between  Sobolev or Besov type norms associated with
free and perturbed evolution flow, see \cite{CasDan-2016}, \cite{DanPie-2005}, \cite{GeoVis-2003}, \cite{IMT-Besov}, \cite{IMT-bilinear}) .
A typical gradient estimate for the classical heat equation
$$ \partial_t u(t,x) -\Delta u(t,x) = 0, \quad  t\in(0,\infty),\quad x\in \R^n$$
is the following one
\begin{equation}
\label{EQ:fe}
\|\nabla u(t)\|_{L^p(\R^n)}
\le
C t^{-\frac{1}{2}}
\|f\|_{L^p(\R^n)}
\end{equation}
valid for any $t>0$ and $1 \le p \le \infty.$
The proof follows immediately from the explicit representation formula of $u.$
The estimate \eqref{EQ:fe} becomes
\begin{equation}
\label{EQ:p-grad}
\|\nabla u(t)\|_{L^p(\Omega)}
\le
C t^{-\frac{1}{2}}
\|f\|_{L^p(\Omega)}
\end{equation}
in
the case of initial boundary value problem with Dirichlet boundary condition in domain $\Omega$.
This estimate is true for any $t>0$ and $1\le p \le \infty,$ when $\Omega$ is a half space.
Again key point in the proof of \eqref{EQ:p-grad} is  the explicit representation formula of $u$.

However,  the question whether an optimal  gradient estimate similar to \eqref{EQ:p-grad} is true for the linear heat flow in arbitrary exterior domain with Dirichlet boundary condition for any $t >0$ and $1 \le p \le \infty$ seems to remain without complete answer due to our knowledge.

Surprisingly, more information and in particular answers to this question can be found in  the case of Stokes equations in exterior domain.
In fact, the estimate \eqref{EQ:p-grad} has appropriate modification in the case of initial boundary value problem for the Stokes equation in  exterior domain $\Omega \subset \mathbb{R}^n.$
The case of Stokes equations with Dirichlet boundary condition
\begin{equation*}
\begin{cases}
	\partial_t u(t,x) -\Delta u(t,x) + \nabla p = 0, \ \ {\rm div} u =0, \quad & t\in(0,\infty),\quad x\in\Omega,\\
	u(t,x) = 0, & t\in(0,\infty),\quad x\in\partial\Omega,\\
	u(0,x) = f(x), & x\in \Omega
\end{cases}
\end{equation*}
is studied in \cite{MarSol-1997}, \cite{DSh99}, where the gradient estimate

\begin{equation}
\label{EQ:fe1b}
\|\nabla u(t)\|_{L^p(\Omega)}
\le
\begin{cases}
C t^{-\frac{1}{2}} \|f\|_{L^p(\Omega)} \quad &\text{for }0 < t \le 1,\\
C t^{-\mu}
\|f\|_{L^p(\Omega)} &\text{for }t\ge 1
\end{cases}
\end{equation}
with
\begin{equation*}
\mu=
\begin{cases}
\frac{1}{2} \quad &\text{if $1 \le p \le n$,}\\
\frac{n}{2p} \quad &\text{if $n < p \le \infty$}
\end{cases}
\end{equation*}
is verified. The optimality of the estimate of \eqref{EQ:fe1b} is discussed in \cite{MarSol-1997}, where the authors show that estimate of type
 \begin{equation}
 \label{EQ:opti}
\|\nabla u(t)\|_{L^p(\Omega)}
\le
C t^{-\mu-\delta}
\|f\|_{L^p(\Omega)},
\quad t >1, \quad \delta >0
\end{equation}
is not true.

Estimate \eqref{EQ:fe1b} for the
Dirichlet problem of heat equation in $\Omega$:
\begin{equation}
\label{EQ:DP}
\begin{cases}
	\partial_t u(t,x) -\Delta u(t,x) = 0, \quad & t\in(0,\infty),\quad x\in\Omega,\\
	u(t,x) = 0, & t\in(0,\infty),\quad x\in\partial\Omega,\\
	u(0,x) = f(x), & x\in \Omega
\end{cases}
\end{equation}
is studied for several situations.
The case  $1 \le p \leq 2$
is studied in \cite{IMT-bdd}, \cite{Ouh_2005}, where the estimate \eqref{EQ:p-grad} is proved for any $t>0$ in an arbitrary open set.
On the other hand, the situation in the case $p>2$ is more complicated.
The case of the Ornstein - Uhlenbeck semigroup, including heat semigroup as a special case, is considered for $t>0$ and $1 < p < \infty$ in \cite{GHHW05}.
The case of parabolic equation is considered in
\cite{Lun_1995}, \cite{ForMetPri-2004} and 
 the results obtained in these works imply that 
the bounded classical solution
satisfies the gradient estimate \eqref{EQ:p-grad} with $p=\infty$ for $0<t\le T$.

Our first goal shall be to prove  the gradient estimate \eqref{EQ:fe1b}  for the heat flow in general exterior domain with Dirichlet boundary condition.
Our second main point is suggested by the following observation. One can expect that the estimate \eqref{EQ:p-grad} shall be true at least in  the special case, when  $\Omega$ is an exterior of a ball, since one can have an explicit representation of the solution $u$. Note that the estimate \eqref{EQ:p-grad} is stronger than the estimate \eqref{EQ:fe1b} for $p > n$ as $t \to \infty.$
For this, our next step shall be to show that  \eqref{EQ:p-grad} is not fulfilled  when $\Omega$ is an exterior of a ball. In this case,
denoting by $u(t;f)$ the heat flow solution to \eqref{EQ:DP} with initial data $f$,
we can show that
$$ 0 < \sup_{t>0, f \in L^p(\Omega), \|f\|_{L^p(\Omega)} = 1}
t^\mu \|\nabla u(t;f)\|_{L^p(\Omega)} < \infty $$ for any $1 \le p \le \infty.$
The right inequality follows from the gradient estimate \eqref{EQ:fe1b}. The left inequality, i.e. the positivity of the supremum gives variational characterization of the best constant in  \eqref{EQ:fe1b} and implies the optimality of the gradient estimate at least when $\Omega$ is the exterior of a ball.

\section{Assumptions and main results}
Let $\Omega$ be the exterior domain of a compact set,
whose boundary is assumed to be sufficiently smooth.

As we have promised in the introduction, our first result concerns gradient estimate for the heat flow outside $\Omega$ with Dirichlet boundary condition.

\begin{thm}
\label{thm:grad}
Let $n\ge2$ and $\Omega$ be the exterior domain in $\mathbb R^n$ of a compact set with $C^{1,1}$ boundary.
Then,
for any $1 \le p < \infty$, there exists a constant $C>0$ such that
the solution $u$ of \eqref{EQ:DP} satisfies
\begin{equation}
\label{EQ:grad-heat}
\|\nabla u(t)\|_{L^p(\Omega)}
\le
\begin{cases}
C t^{-\frac{1}{2}}
\|f\|_{L^p(\Omega)}\quad &\text{for $0<t \le 1$,}\\
C t^{-\mu}
\|f\|_{L^p(\Omega)}\quad &\text{for $t>1$}
\end{cases}
\end{equation}
for any $f\in L^p(\Omega)$, where the exponent $\mu$ is given by
\begin{equation}\label{EQ:mu}
\mu=
\begin{cases}
\frac{1}{2} \quad &\text{if $1 \le p \le n$,}\\
\frac{n}{2p} \quad &\text{if $n < p < \infty$}.
\end{cases}
\end{equation}
\end{thm}

\begin{rem} Since we consider boundary with weak regularity, it is not clear whether the gradient estimate \eqref{EQ:grad-heat} with $p=\infty$ is true
for any $f\in L^\infty(\Omega)$ due to our knowledge.
However the gradient estimate is true for classical solutions.
In fact, the bounded classical solutions to Dirichlet problem of parabolic equations in bounded or unbounded domains with sufficiently smooth boundary satisfy the estimate
\[
\|\nabla u(t)\|_{L^\infty(\Omega)}
\le C_T t^{-\frac{1}{2}} \|f\|_{L^\infty(\Omega)},\quad 0<t\le T
\]
(see \cite{Lun_1995}, \cite{ForMetPri-2004}, and references therein).

In particular,
it can be proved by using the above gradient estimate combined with $L^p$-$L^q$-estimates of Lemma \ref{lem:LpLq} below
imply that the bounded classical solutions of heat equation \eqref{EQ:DP} satisfy the estimate
\[
\|\nabla u(t)\|_{L^\infty(\Omega)}
\le C  \|f\|_{L^\infty(\Omega)},\quad t>1.
\]
\end{rem}

\begin{rem}
Assuming we have  $C^{1,1}$ boundary, the authors in \cite{GHHW05} have established that there exists $\omega \geq 0$ so that the following estimate
\begin{align}\label{eq.loc1}
   \|\nabla u(t)\|_{L^p(\Omega)}
\le C t^{-\frac{1}{2}} e^{\omega t}
\|f\|_{L^p(\Omega)},\quad t>0,\quad 1 < p < \infty
\end{align}
holds. From this estimate we can see that
\begin{align*}
   \|\nabla u(t)\|_{L^p(\Omega)}
\le C t^{-\frac{1}{2}} \|f\|_{L^p(\Omega)},\quad 0 < t \leq 1, \quad  1 < p < \infty
\end{align*}
is fulfilled, i.e. we have \eqref{EQ:grad-heat} for small values of $t$ and $1 < p < \infty.$ The estimate \eqref{EQ:grad-heat} for small values of $t$  and for the endpoint  case $p=1$ can be deduced from the results in \cite{IMT-bdd}, \cite{Ouh_2005}.
\end{rem}

\begin{rem}
In the case of Neumann boundary condition, the estimate \eqref{EQ:grad-heat} can be replaced by its stronger version \eqref{EQ:p-grad}, i.e. we have the same estimate as in the case of the whole space $\mathbb{R}^n.$
Indeed, in this case,
the gradient estimate \eqref{EQ:p-grad} holds for any $t>0$ and $1\le p \le \infty$ (see, e.g., \cite{Ish-2009} and \cite{Tan-arxiv}).
\end{rem}

Next, we discuss the optimality of time decay rates in the gradient estimates \eqref{EQ:grad-heat} as $t\to \infty$.

\begin{defi}\label{defi:optimal}
We say that the gradient estimate \eqref{EQ:grad-heat} is optimal if
there exist sequences $\{f_m\}_{m\in\mathbb N}\subset L^p(\Omega)$ and $\{t_m\}_{m\in\mathbb N}$ such that
\[
t_m > 0 \quad \text{for }m\in\mathbb N,\quad t_m \to \infty \quad \text{as }m\to\infty
\]
and
\begin{equation*}
\limsup_{m\to \infty} \frac{t_m^\mu \|\nabla u_m(t_m)\|_{L^p(\Omega)}}{\|f_m\|_{L^p(\Omega)}}
>0,
\end{equation*}
where $u_m$ is a solution to \eqref{EQ:DP} with initial data $f_m$ and the exponent $\mu$ is given by \eqref{EQ:mu}.
\end{defi}

\begin{rem}
If we can verify the optimality of \eqref{EQ:grad-heat} in the sense of Definition \ref{defi:optimal},
then we can assert
$$ \sup_{t>0, f \in L^p(\Omega), \|f\|_{L^p(\Omega)} = 1}
t^\mu \|\nabla u(t;f)\|_{L^p(\Omega)}  $$ is a well-defined positive number that gives a variational characterization of the best constant $C=C(\Omega,p)$ in \eqref{EQ:grad-heat}.
\end{rem}

Our result on the optimality is the following.
To simplify the proof, we shall fix the space dimension $n=3$.

\begin{thm}\label{thm:optimal}
Let $n=3$ and $\Omega$ be the exterior domain of a ball. Then, for any $1\le p \le \infty$, the gradient estimate \eqref{EQ:grad-heat} is optimal in
the sense of Definition \ref{defi:optimal}.
\end{thm}

\begin{rem}
Note that the optimality of estimate of type \eqref{EQ:grad-heat} is verified in the context discussed in   \cite{MarSol-1997} for arbitrary exterior domains $\Omega$. However, optimality treated in this work means that

$$ \limsup_{t \to \infty} \sup_{ f \in L^p(\Omega), \|f\|_{L^p(\Omega)} = 1} t^{\mu+ \delta} \|\nabla u(t;f)\|_{L^p(\Omega)} =+\infty$$
for any $\delta>0.$
In other words, \eqref{EQ:opti} is not true for any $\delta>0$. This optimality is weaker than that of Definition \ref{defi:optimal}.
\end{rem}

\begin{rem}
In the case $p=\infty$, when $n\ge 3$ and $\Omega$ is the exterior domain of a compact connected set with $C^{1,1}$ boundary, we can obtain the estimate
\begin{equation}\label{EQ:opti-infi}
\|\nabla e^{-tH}\|_{L^\infty(\Omega)\to L^\infty(\Omega)} \ge C
\end{equation}
for any $t>1$, where $e^{-tH}$ is the semigroup generated by the Dirichlet Laplacian $H=-\Delta$ on $\Omega$ (see Appendix \ref{App:A}).
This is stronger than the optimality of Definition \ref{defi:optimal},
since \eqref{EQ:opti-infi} implies the optimality for $p=\infty$ in the sense of Definition \ref{defi:optimal}.
\end{rem}

\begin{rem}
If one compare the optimality result in Theorem \ref{thm:optimal} and the estimate \eqref{eq.loc1}, then one can deduce that $\omega >0$ in \eqref{eq.loc1}, therefore the assertion (b) of Theorem 3.1 in \cite{GHHW05} holds for some $\omega >0$ and it is not true for $\omega =0.$
\end{rem}

The behavior of the derivatives of solutions to heat equations is
not only of interest itself, but also has some applications.
We can present an a priori estimate concerning the bilinear estimates for  Dirichlet Laplacian $H$  on $L^2(\Omega)$ as an application of gradient estimates obtained in the previous theorems.

The bilinear estimates are of great importance
in the study of the Dirichlet problem for nonlinear partial differential equations.
Recently, it is revealed in \cite{IMT-bilinear} that
the gradient estimate \eqref{EQ:p-grad} implies bilinear estimates in Besov spaces.
According to the result, by proving the gradient estimates in exterior domains,
we obtain the bilinear estimates in exterior domains.

To state this result,
we recall the definition of the homogeneous Besov spaces  $\dot{B}^s_{p,q}(H)$ with norm
\[
\|f\|_{\dot{B}^s_{p,q}(H)}
=
\left\{
\sum_{j=-\infty}^{\infty}
\left(2^{sj} \|\phi_j(\sqrt{H}) f\|_{L^p(\Omega)}\right)^q
\right\}^{\frac{1}{q}},
\]
where   $s\in\mathbb R,$ $1 \le p,q\le \infty$, and we have used the functional calculus for the self - adjoint operator $H$ in combination with Paley - Littlewood partition
in the right side of this identity.

\begin{cor}\label{cor:bilinear}
Let  $\Omega$ be the exterior domain in $\mathbb R^n$ of a compact set with $C^{1,1}$ boundary.
Let $0<s<2$ and $p$, $p_1,p_2,p_3,p_4$ and
$q$ be such that
$$
1\le p, p_1,p_2,p_3,p_4\le n, \quad
1\le q\le \infty
\quad \text{and}
\quad \frac{1}{p}=\frac{1}{p_1}+\frac{1}{p_2}=\frac{1}{p_3}+\frac{1}{p_4}.
$$
Then there exists a constant $C>0$ such that
\begin{equation*}
\|fg\|_{\dot{B}^s_{p,q}(H)}
\le
C\left(
\|f\|_{\dot{B}^s_{p_1,q}(H)}
\|g\|_{L^{p_2}(\Omega)}
+
\|f\|_{L^{p_3}(\Omega)}
\|g\|_{\dot{B}^s_{p_4,q}(H)}
\right)
\end{equation*}
for any
$f\in \dot{B}^s_{p_1, q}(H)\cap L^{p_3}(\Omega)$ and $g\in \dot{B}^s_{p_4, q}(H)\cap
L^{p_2}(\Omega)$.
\end{cor}

\section{Key estimates}

In this section we prepare key estimates for solutions of heat equations \eqref{EQ:DP}.
The first one is the result on $L^p$-$L^q$-estimates.

\begin{lem}\label{lem:LpLq}
Let $\Omega$ be an open set in $\mathbb R^n$ and $1 \le p \le q \le \infty$. Then
there exists a constant $C>0$ such that
\begin{equation*}
\|u(t)\|_{L^q(\Omega)} \le C t^{-\frac{n}{2}(\frac{1}{p}-\frac{1}{q})} \|f \|_{L^p(\Omega)}
\end{equation*}
for any $t>0$ and $f \in L^p(\Omega)$.
\end{lem}

For the proof, we refer to Proposition 3.1 in \cite{IMT-bdd} (see also Section 6.3 in \cite{Ouh_2005}). \\

The second one is the result on the gradient estimates for $1\le p\le 2$.

\begin{lem}\label{lem:grad}
Let $\Omega$ be an open set in $\mathbb R^n$ and $1 \le p \le 2$. Then
\[
\|\nabla u(t)\|_{L^p(\Omega)}
\le
C t^{-\frac{1}{2}}
\|f\|_{L^p(\Omega)}
\]
for any $t>0$ and $f \in L^p(\Omega)$.
\end{lem}

For the proof, we refer to Theorem 1.2 in \cite{IMT-bdd} (see also Theorem 6.19 in \cite{Ouh_2005}).


\section{Proof of Theorem \ref{thm:grad}}
In this section we prove Theorem \ref{thm:grad}.
The case $p=1$ is proved in Lemma \ref{lem:grad}.
Hence, in order to obtain \eqref{EQ:grad-heat} for any $1\le p < \infty$, it suffices to prove the case $n \le p < \infty$ by density and interpolation argument: For any $n \le p<\infty$, there exists a constant $C>0$ such that
\begin{equation}
\label{EQ:aim}
\|\nabla u(t)\|_{L^p(\Omega)}
\le
\begin{cases}
C t^{-\frac{1}{2}}\|f\|_{L^p(\Omega)}\quad &\text{for $0<t\le 1$},\\
C t^{-\frac{n}{2p}}\|f\|_{L^p(\Omega)} &\text{for $t > 1$}
\end{cases}
\end{equation}
for any $f \in C^\infty_0(\Omega)$.
Let us choose $L>0$ such that
\begin{equation}\label{EQ:L}
\mathbb R^n \setminus \Omega \subset \{|x| < L\}.
\end{equation}
Putting
\[
\Omega_{L+2}:=\Omega \cap \{|x|< L+2\},
\]
we estimate
\begin{equation}\label{EQ:1}
\|\nabla u(t)\|_{L^p(\Omega)}
\le \|\nabla u(t)\|_{L^p(\Omega_{L+2})}
+
\|\nabla u(t)\|_{L^p(\{|x|\ge L+2\})}.
\end{equation}
As to the first term,
we can obtain
\begin{equation}
\label{EQ:first}
\|\nabla u(t)\|_{L^p(\Omega_{L+2})}
\le
\begin{cases}
C t^{-\frac{1}{2}}\|f\|_{L^p(\Omega)}\quad &\text{for $0<t\le 1$},\\
C t^{-\frac{n}{2p}}\|f\|_{L^p(\Omega)} &\text{for $t > 1$}
\end{cases}
\end{equation}
by using Lemmas \ref{lem:GN} and \ref{lem:B} in Appendix B.
In fact, noting that
\[
u(t) \in W^{2,p}(\Omega)\cap W^{1,p}_0(\Omega)
\]
for any $t>0$ and $f \in C^\infty_0(\Omega)$, we can apply Lemmas \ref{lem:GN} and \ref{lem:B}
to estimate
\begin{equation}\label{EQ:11}
\begin{split}
\|\nabla u(t)\|_{L^p(\Omega_{L+2})}
& \le
C_1
\|D^\alpha u(t)\|_{L^p(\Omega_{L+2})}^{\frac{1}{2}} \|u(t)\|_{L^p(\Omega_{L+2})}^{\frac{1}{2}}
+ C_2 \|u(t)\|_{L^p(\Omega_{L+2})}\\
& \le
C \left(
\|\Delta u(t)\|_{L^p(\Omega)}^{\frac{1}{2}} \|u(t)\|_{L^p(\Omega)}^{\frac{1}{2}}
+  \|u(t)\|_{L^p(\Omega_{L+2})}
\right),
\end{split}
\end{equation}
where $\alpha$ is a multi-index with $|\alpha|=2$.
Since
\[
\|\Delta u(t)\|_{L^p(\Omega)}^{\frac{1}{2}} \|u(t)\|_{L^p(\Omega)}^{\frac{1}{2}}
\le C t^{-\frac{1}{2}} \|f\|_{L^p(\Omega)}
\]
and
\[
\|u(t)\|_{L^p(\Omega_{L+2})}
\le
C \|u(t)\|_{L^\infty(\Omega)} \le
C t^{-\frac{n}{2p}} \|f\|_{L^p(\Omega)}
\]
for any $t>0$ by Lemma \ref{lem:LpLq},
the right hand side in \eqref{EQ:11} is estimated as
\[
\|\Delta u(t)\|_{L^p(\Omega)}^{\frac{1}{2}} \|u(t)\|_{L^p(\Omega)}^{\frac{1}{2}}
+  \|u(t)\|_{L^p(\Omega_{L+2})}
\le C \max (t^{-\frac{1}{2}}, t^{-\frac{n}{2p}}) \|f\|_{L^p(\Omega)}
\]
for any $t>0$. Therefore we obtain the required estimates \eqref{EQ:first}.
Thus all we have to do is to estimate the second term in \eqref{EQ:1} as follows:
\begin{equation}
\label{EQ:second}
\|\nabla u(t)\|_{L^p(\{|x|> L+2\})}
\le
\begin{cases}
C t^{-\frac{1}{2}}\|f\|_{L^p(\Omega)}\quad &\text{for $0<t\le 1$},\\
C t^{-\frac{n}{2p}}\|f\|_{L^p(\Omega)} &\text{for $t > 1$}.
\end{cases}
\end{equation}
We divide the proof of \eqref{EQ:second} into two cases:
$0<t \le 1$ and $t>1$. \\

\noindent
{\bf The case $0<t \le 1$.}
We denote by $\chi_L$ a smooth function on $\mathbb R^n$ such that
\begin{equation}\label{EQ:chi_L}
\chi_L(x)=
\begin{cases}
1\quad \text{for $|x|\ge L+1$},\\
0\quad \text{for $|x|\le L$},
\end{cases}
\end{equation}
and have
\[
u(t,x)= \chi_L(x) u(t,x),\quad |x|\ge L+2.
\]
Let us decompose $\chi_L u(t)$ into
\begin{equation}\label{EQ:decom}
\chi_L u(t) = v_1(t)-v_2(t)
\end{equation}
for $0<t \le 1$. Here
$v_1(t)$ is the solution to the Cauchy problem of heat equation in $\mathbb R^n$:
\begin{equation*}
\begin{cases}
	\partial_t v_1(t,x) -\Delta v_1(t,x) = 0, \quad & t\in(0,1],\quad x\in\mathbb R^n,\\
	v_1(0,x) = \chi_L(x)f(x), & x\in \mathbb R^n,
\end{cases}
\end{equation*}
and $v_2(t)$ is the solution to the Cauchy problem of heat equation in $\mathbb R^n$:
\begin{equation*}
\begin{cases}
	\partial_t v_2(t,x) -\Delta v_2(t,x) = F(t,x), \quad & t\in(0,1],\quad x\in\mathbb R^n,\\
	v_2(0,x) = 0, & x\in \mathbb R^n,
\end{cases}
\end{equation*}
where
\begin{equation}\label{EQ:F}
F(t,x)=-2\nabla \chi_L(x) \cdot \nabla u(t,x) + (\Delta \chi_L(x)) u(t,x).
\end{equation}
It is easily proved that
\begin{equation}\label{EQ:v_1-grad}
\|\nabla v_1(t)\|_{L^p(\{|x|> L+2\})} \le C t^{-\frac{1}{2}}
\|f\|_{L^p(\Omega)}
\end{equation}
for any $0<t\le1$.
Hence it is sufficient to show that
\begin{equation}\label{EQ:v_2-grad}
\|\nabla v_2(t)\|_{L^p(\{|x|> L+2\})} \le C t^{-\frac{1}{2}}
\|f\|_{L^p(\Omega)}
\end{equation}
for any $0<t\le1$.
Letting $e^{t\Delta}$ be the semigroup generated by $-\Delta$ on $\mathbb R^n$,
we write $v_2(t)$ as
\[
v_2(t,x) = \int_0^t e^{(t-s)\Delta}F(s,x)\,ds
\]
for $0<t\le1$ and $x\in \mathbb R^n$.
Denoting by
$G(t,x-y)$ the kernel of $e^{t\Delta}$, we estimate
\begin{equation}\label{EQ:v-1}
\begin{split}
& \|\nabla v_2(t)\|_{L^p(\{|x|> L+2\})}\\
\le \, &
\int_0^t \|\nabla e^{(t-s)\Delta}F(s,\cdot)\|_{L^p(\{|x|> L+2\}}\,ds\\
\le \, &
\int_0^t \int_{L<|y|\le L+1}\|\nabla_x G(t-s,x-y)\|_{L^p(\{|x|> L+2\})} |F(s,y)|\,dy\,ds.
\end{split}
\end{equation}
Here we note that
\begin{equation*}
\begin{split}
\left|\nabla_x G(t,x-y)\right|
& \le C t^{-\frac{n+1}{2}}
\frac{|x-y|}{2\sqrt{t}} e^{-\frac{|x-y|^2}{4t}}\\
& \le C t^{-\frac{n+1}{2}} \left(1 + \frac{|x-y|^2}{t}\right)^{-\frac{n+1}{2}}\\
& = C(t+|x-y|^2)^{-\frac{n+1}{2}}
\end{split}
\end{equation*}
for any $t>0$ and $x\in \mathbb R^n$.
In particular, if $|x|\ge L+2$ and $|y| \le L+1$, then
\[
|x-y| \ge |x| - |y|
\ge |x| - \frac{L+1}{L+2} |x|
= \frac{1}{L+2}|x|,
\]
and hence,
\begin{equation*} 
|\nabla_x G(t-s,x-y)| \le C \{(t-s)+|x|^2\}^{-\frac{n+1}{2}}
\end{equation*}
for any $0<s<t$. Therefore we deduce that
\begin{equation}\label{EQ:v-2}
\|\nabla_x G(t-s,x-y)\|_{L^p(\{|x|> L+2\})}
\le C\{1 + (t-s)\}^{-\frac{n+1}{2}+\frac{n}{2p}}
\end{equation}
for any $0<s<t$ and $L<|y|\le L+1$.
Combining \eqref{EQ:v-1} and \eqref{EQ:v-2}, we obtain
\begin{equation}\label{EQ:v-3}
\begin{split}
\|\nabla v_2(t)\|_{L^p(\{|x|> L+2\})}
& \le
\int_0^t \int_{L<|y|\le L+1}\{1 + (t-s)\}^{-\frac{n+1}{2}+\frac{n}{2p}} |F(s,y)|\,dy\,ds\\
& =
\int_0^t \{1 + (t-s)\}^{-\frac{n+1}{2}+\frac{n}{2p}} \|F(s,\cdot)\|_{L^1(\{L<|y|\le L+1\})}\,ds.
\end{split}
\end{equation}
Recalling the definition \eqref{EQ:F} of $F(s,x)$, and using \eqref{EQ:first} and Lemma \ref{lem:LpLq}, we estimate
\begin{equation*} 
\begin{split}
\|F(s,\cdot)\|_{L^1(\{L<|y|\le L+1\})}
& \le C \left(
\|\nabla u(s)\|_{L^1(\{L<|y|\le L+1\})}+\|u(s)\|_{L^1(\{L<|y|\le L+1\})}
\right)\\
&\le C \left(
\|\nabla u(s)\|_{L^p(\{L<|y|\le L+1\})}+\|u(s)\|_{L^p(\{L<|y|\le L+1\})}
\right)\\
& \le C 
s^{-\frac{1}{2}}
\|f\|_{L^p(\Omega)}
\end{split}
\end{equation*}
for any $0<s\le1$.
Combining the above two estimates, we deduce that
\begin{equation*}
\begin{split}
\|\nabla v_2(t)\|_{L^p(\{|x|> L+2\})}
& \le C
\int_0^t \{1 + (t-s)\}^{-\frac{n+1}{2}+\frac{n}{2p}} s^{-\frac{1}{2}}\,ds\cdot \|f\|_{L^p(\Omega)}\\
& \le C
\int_0^t s^{-\frac{1}{2}}\,ds\cdot \|f\|_{L^p(\Omega)}\\
& \le C \|f\|_{L^p(\Omega)}
\end{split}
\end{equation*}
for any $0<t \le 1$, which proves \eqref{EQ:v_2-grad}.
Therefore the estimate \eqref{EQ:second} for any $0<t \le 1$ is proved by \eqref{EQ:v_1-grad} and \eqref{EQ:v_2-grad}. \\

\noindent
{\bf The case $t>1$.}
In a similar way to \eqref{EQ:decom} in the previous case,
we decompose $\chi_L u(t)$ into
\begin{equation*} 
\chi_L u(t) = w_1(t)-w_2(t)
\end{equation*}
for $t\ge1$. Here
$w_1(t)$ is the solution to the Cauchy problem of heat equation in $\mathbb R^n$:
\begin{equation*}
\begin{cases}
	\partial_t w_1(t,x) -\Delta w_1(t,x) = 0, \quad & t\in(1,\infty),\quad x\in\mathbb R^n,\\
	w_1(1,x) = \chi_L(x)u(1,x), & x\in \mathbb R^n,
\end{cases}
\end{equation*}
and $w_2(t)$ is the solution to the Cauchy problem of heat equation in $\mathbb R^n$:
\begin{equation*}
\begin{cases}
	\partial_t w_2(t,x) -\Delta w_2(t,x) = F(t,x), \quad & t\in(1,\infty),\quad x\in\mathbb R^n,\\
	w_2(1,x) = 0, & x\in \mathbb R^n,
\end{cases}
\end{equation*}
where we recall \eqref{EQ:chi_L} and \eqref{EQ:F}.
It is easily proved that
\begin{equation}\label{EQ:w_1-grad}
\|\nabla w_1(t)\|_{L^p(\{|x|> L+2\})} \le C t^{-\frac{1}{2}}
\|f\|_{L^p(\Omega)}
\end{equation}
for any $t>1$.
Hence it is sufficient to show that
\begin{equation}\label{EQ:w_2-grad}
\|\nabla w_2(t)\|_{L^p(\{|x|> L+2\})} \le C t^{-\frac{n}{2p}}
\|f\|_{L^p(\Omega)}
\end{equation}
for any $t>1$.
Writing $w_2(t)$ as
\[
w_2(t,x) = \int_1^t e^{(t-s)\Delta}F(s,x)\,ds
\]
for $t>1$ and $x\in \mathbb R^n$, we estimate, in a similar way to \eqref{EQ:v-3},
\begin{equation*}
\begin{split}
\|\nabla w_2(t)\|_{L^p(\{|x|> L+2\})}
\le C
\int_1^t \{1 + (t-s)\}^{-\frac{n+1}{2}+\frac{n}{2p}} \|F(s,\cdot)\|_{L^1(\{L<|y|\le L+1\})}\,ds.
\end{split}
\end{equation*}
Recalling the definition \eqref{EQ:F} of $F(s,x)$, and using \eqref{EQ:first} and Lemma \ref{lem:LpLq}, we estimate
\begin{equation*} 
\begin{split}
\|F(s,\cdot)\|_{L^1(\{L<|y|\le L+1\})}
& \le C \left(
\|\nabla u(s)\|_{L^1(\{L<|y|\le L+1\})}+\|u(s)\|_{L^1(\{L<|y|\le L+1\})}
\right)\\
&\le C \left(
\|\nabla u(s)\|_{L^p(\{L<|y|\le L+1\})}+\|u(s)\|_{L^\infty(\{L<|y|\le L+1\})}
\right)\\
& \le C 
s^{-\frac{n}{2p}}
\|f\|_{L^p(\Omega)}
\end{split}
\end{equation*}
for any $s>1$.
Combining the above two estimates, we deduce that
\begin{equation*}
\begin{split}
\|\nabla w_2(t)\|_{L^p(\{|x|> L+2\})}
 \le C
\int_1^t \{1 + (t-s)\}^{-\frac{n+1}{2}+\frac{n}{2p}} s^{-\frac{n}{2p}} \,ds\cdot \|f\|_{L^p(\Omega)}.
\end{split}
\end{equation*}
for any $t>1$.
For $1 < t < 2$ we use the inequality
\begin{equation*} 
   \int_1^{t} \{1 + (t-s)\}^{-\frac{n+1}{2}+\frac{n}{2p}} s^{-\frac{n}{2p}}\,ds \leq  \int_1^{2}s^{-\frac{n}{2p}}\,ds \leq Ct^{-\frac{1}{2}}.
\end{equation*}
For $t>2$ and $p \geq n$, we have
\begin{equation*}
\begin{split}
\int_1^{\frac{t}{2}} \{1 + (t-s)\}^{-\frac{n+1}{2}+\frac{n}{2p}} s^{-\frac{n}{2p}}\,ds
& \le C t^{-\frac{n+1}{2}+\frac{n}{2p}}\int_1^{\frac{t}{2}} s^{-\frac{n}{2p}}\,ds\\
& \le C t^{-\frac{n}{2}+\frac{1}{2}} \le Ct^{-\frac{1}{2}}
\end{split}
\end{equation*}
and
\begin{equation*}
\begin{split}
\int_{\frac{t}{2}}^t \{1 + (t-s)\}^{-\frac{n+1}{2}+\frac{n}{2p}} s^{-\frac{n}{2p}}\,ds
& \le C t^{-\frac{n}{2p}} \int_{\frac{t}{2}}^t \{1 + (t-s)\}^{-\frac{n+1}{2}+\frac{n}{2p}} \,ds\\
& \le C t^{-\frac{n}{2p}} \left(1 + t^{-\frac{n-1}{2} + \frac{n}{2p}} \right)\\
& \le C \left(t^{-\frac{n}{2p}} + t^{-\frac{n-1}{2}} \right)
\le C t^{-\frac{n}{2p}}.
\end{split}
\end{equation*}
Hence
we obtain the estimate \eqref{EQ:w_2-grad}
for any $t>1$.
Therefore the estimate \eqref{EQ:second} for any $t>1$ is proved by \eqref{EQ:w_1-grad} and  \eqref{EQ:w_2-grad}.

Thus, combining \eqref{EQ:1} with \eqref{EQ:first} and \eqref{EQ:second},  we conclude the estimates \eqref{EQ:aim}. The proof of Theorem \ref{thm:grad} is complete.

\section{Proof of Theorem \ref{thm:optimal}}
In this section we discuss the optimality of time decay rates of estimates \eqref{EQ:grad-heat} in Theorem \ref{thm:grad}. Let $\Omega$ be the exterior domain in $\mathbb R^3$ determined by
\[
\Omega = \{x\in\mathbb R^3:|x|>1\}.
\]
We prove the optimality in the sense of Definition \ref{defi:optimal} for heat equation \eqref{EQ:DP} with a radial initial data $f$ on $\Omega$.

\begin{proof}
Let $f$ be a radial function on $\Omega$.
Since $u(t)$ is also radial,
we write
\[
F(r):=f(x),\quad U(t,r):=u(t,x)
\]
for $t>0$ and $r=|x|$.
We rewrite the problem \eqref{EQ:DP} to the following problem by the polar coordinates and making change $v(t,r)=(r+1)U(t,r+1)$:
\begin{equation}
\label{EQ:DP-rad}
\begin{cases}
	\partial_t v(t,r) - \partial^2_r v(t,r)
	= 0, \quad & t\in(0,\infty),\quad r\in (0, \infty),\\
	v(t,r) = 0, & t\in(0,\infty),\quad r=0,\\
	v(0,r) = g(r), & r\in (0, \infty),
\end{cases}
\end{equation}
where $g(r)=(r+1) F(r+1)$ and $r=|x|$.
Then
solutions $v$ to \eqref{EQ:DP-rad} and the derivative $\partial_r v$ can be
represented as
\begin{equation}
\label{EQ:repre}
\begin{split}
v(t,r)
 =
(4 \pi t)^{-\frac{1}{2}}
\int^\infty_0
\left\{
e^{-\frac{(r-s)^2}{4t}}
- e^{-\frac{(r+s)^2}{4t}} \right\}g(s)\, ds,
\end{split}
\end{equation}
\begin{equation}
\label{EQ:repre-d}
\begin{split}
\partial _r v(t,r)
& =
(4 \pi t)^{-\frac{1}{2}}
\int^\infty_0
\left\{
- \frac{r-s}{2t} e^{-\frac{(r-s)^2}{4t}}
+ \frac{r+s}{2t} e^{-\frac{(r+s)^2}{4t}} \right\}g(s)\, ds
\end{split}
\end{equation}
for $t>0$ and $r>0$.
Furthermore, noting that $u(t,x)=U(t,r) = r^{-1}v(t, r-1)$,
we write
\begin{equation}
\label{EQ:rewrite}
\begin{split}
&\|\nabla u(t) \|_{L^p(\Omega)}\\
= &
(4\pi)^{\frac{1}{p}}
\left(
\int^\infty_1
\left|\partial_r U(t,r) \right|^p r^{2}\, dr
\right)^{\frac{1}{p}}\\
= &
(4\pi)^{\frac{1}{p}}
\left(
\int^\infty_1
\left|
\partial_r \left(r^{-1}v(t,r-1)\right)
\right|^p r^{2}\, dr
\right)^{\frac{1}{p}}\\
= &
(4\pi)^{\frac{1}{p}}
\left(
\int^\infty_1
\left| -r^{-2} v(t,r-1) +
r^{-1}\partial_r v(t,r-1)
\right|^p r^{2}\, dr
\right)^{\frac{1}{p}}\\
= &
(4\pi)^{\frac{1}{p}}
\left(
\int^\infty_0
\left| -(r+1)^{-2+\frac{2}{p}} v(t,r) +
(r+1)^{-1+\frac{2}{p}}\partial_r v(t,r)
\right|^p \, dr
\right)^{\frac{1}{p}}\\
= &
(4\pi)^{\frac{1}{p}}
\left\| -(r+1)^{-2+\frac{2}{p}} v(t) + (r+1)^{-1+\frac{2}{p}} \partial_r v(t) \right\|_{L^p(0, \infty)}.
\end{split}
\end{equation}
In order to prove the optimality,
we choose appropriate initial data $f_m$ and
estimate from below the  quantity from Definition \ref{defi:optimal}:
\begin{equation*}
\frac{t_m^{\mu} \|\nabla u_m(t_m)\|_{L^p(\Omega)}}{\|f_m\|_{L^p(\Omega)}}
\end{equation*}
for $m\in\mathbb N$, where the exponent $\mu$ is defined in \eqref{EQ:mu}.
We divide the proof into two cases: $1\le p\le 3$ and $3<p\le \infty$.\\

\noindent
{\bf The case $1\le p\le 3$.}
We take
\begin{equation*} 
    t_m = m^2
\end{equation*}
for $m\in\mathbb N$,
and define the initial data as follows
\begin{equation}
\label{EQ:f_m}
f_{m}(x) :=
\begin{cases}
C_{m}|x|^{-1},\quad &r\in (m+1,2m+1],\\
0, &\text{otherwise}.\\
\end{cases}
\end{equation}
Here we choose the constant $C_{m}$ such that
\begin{equation}\label{EQ:f_m2}
C_{m} > 0\quad
\text{and} \quad
\|f_{m}\|_{L^p(\Omega)} = 1.
\end{equation}
Then we have
\begin{equation}
\label{EQ:g_m}
g_{m}(r) =
\begin{cases}
C_{m},\quad &r\in (m,2m],\\
0, &\text{otherwise},\\
\end{cases}
\end{equation}
and
\begin{equation}\label{EQ:C_m}
C_{m} \sim m^{1-\frac{3}{p}}
\end{equation}
as $m\to\infty$.
Let us denote by $u_{m}$ and $v_{m}$ the solutions to \eqref{EQ:DP} and \eqref{EQ:DP-rad} with initial data $f_{m}$ and $g_{m}$, respectively.
By the equality \eqref{EQ:rewrite}, we write
\begin{equation*}
\|\nabla u_{m}(t) \|_{L^p(\Omega)}
=
(4\pi)^{\frac{1}{p}}
\left\| -(r+1)^{-2+\frac{2}{p}} v_{m}(t) + (r+1)^{-1+\frac{2}{p}} \partial_r v_{m}(t) \right\|_{L^p(0,\infty)}.
\end{equation*}
Letting $t>0$ and $s>0$ be fixed,
we see that the function
\[
e^{-\frac{(r-s)^2}{4t}}
- e^{-\frac{(r+s)^2}{4t}}, \quad r>0,
\]
is monotonically decreasing with respect to $r \in [\sqrt{2t}+s, \infty)$. Hence, noting from \eqref{EQ:g_m} that $g_m \ge 0$ and $m\le s \le 2m$,
we have
\begin{equation*}
v_{m}(t,r) \ge 0\quad\text{and}\quad \partial_r v_{m}(t,r) \le 0
\end{equation*}
for any $r \in [\sqrt{2t}+2m, \infty)$.
Thanks to this observation, we estimate from below
\[
\|\nabla u_{m}(t)\|_{L^p(\Omega)}
\ge
\left\| (r+1)^{-2+\frac{2}{p}} v_{m}(t)  \right\|_{L^p(\sqrt{2t}+2m,\infty)}.
\]
Taking $t=t_m=m^2$, we write
\begin{equation}
\label{EQ:u_m}
\|\nabla u_{m}(t_m)\|_{L^p(\Omega)}
\ge
\left\| (r+1)^{-2+\frac{2}{p}} v_{m}(m^2)  \right\|_{L^p(c_0m,\infty)},
\end{equation}
where $c_0=2+\sqrt{2}$.
From the representation \eqref{EQ:repre} and definition \eqref{EQ:g_m} of $g_m$, the right hand side is estimated as
\begin{equation}\label{EQ:u_m1}
\begin{split}
& \left\| (r+1)^{-2+\frac{2}{p}} v_{m}(m^2)  \right\|_{L^p(c_0m,\infty)}\\
\ge  & \,
C\cdot C_m m^{-1}
\left\| (r+1)^{-2+\frac{2}{p}}
\int_m^{2m}
\left\{
e^{-\frac{(r-s)^2}{4m^2}}
- e^{-\frac{(r+s)^2}{4m^2}}
\right\}
\,ds  \right\|_{L^p(c_0m,\infty)}.
\end{split}
\end{equation}
Since
\[
e^{-\frac{(r-s)^2}{4m^2}}
- e^{-\frac{(r+s)^2}{4m^2}}
=
e^{-\frac{(r-s)^2}{4m^2}}
(1-e^{-\frac{rs}{m^2}})
\ge (1-e^{-c_0}) e^{-\frac{c_0^2}{4}}e^{-\frac{s^2}{4m^2}}
\]
for any $r>c_0m$ and $m\le s\le 2m$,
the integral in the right hand side in \eqref{EQ:u_m1} is estimated from below as
\begin{equation}\label{EQ:integral}
\int_m^{2m}
\left\{
e^{-\frac{(r-s)^2}{4m^2}}
- e^{-\frac{(r+s)^2}{4m^2}}
\right\}\,ds
\ge
C
\int_m^{2m}
e^{-\frac{s^2}{4m^2}}\,ds
=
Cm \int_1^2 e^{-\frac{s^2}{4}}\,ds.
\end{equation}
Hence, by combining \eqref{EQ:u_m}--\eqref{EQ:integral}, we estimate from below
\begin{equation}\label{EQ:2}
\left\| (r+1)^{-2+\frac{2}{p}} v_{m}(m^2)  \right\|_{L^p(c_0m,\infty)}
\ge C \cdot C_m \left\| (r+1)^{-2+\frac{2}{p}}  \right\|_{L^p(c_0m,2c_0 m)}.
\end{equation}
Hence, noting from \eqref{EQ:C_m} that
\[
C_m \left\| (r+1)^{-2+\frac{2}{p}}  \right\|_{L^p(c_0m,2c_0 m)}
\sim m^{1-\frac{3}{p}} m^{-2+\frac{3}{p}} = m^{-1} = t_m^{-\frac{1}{2}}
\]
as $m\to \infty$, we deduce from \eqref{EQ:u_m}--\eqref{EQ:2} that
\begin{equation}
\label{EQ:u_m3}
\|\nabla u_{m}(t_m)\|_{L^p(\Omega)}
\ge C t_m^{-\frac{1}{2}}
\end{equation}
for sufficiently large $m\in\mathbb N$, where the constant $C>0$ is independent of $m$. %
By combining
\eqref{EQ:f_m2} and \eqref{EQ:u_m3},
we conclude that
\[
\limsup_{m\to\infty}\frac{t_m^{\frac{1}{2}}\|\nabla u_{m}(t_m)\|_{L^p(\Omega)}}{\|f_m\|_{L^p(\Omega)}} >0.
\]
Thus the optimality for $1\le p\le 3$ is proved. \\

\noindent
{\bf The case $3 < p\le \infty$.}
Recalling the equality \eqref{EQ:rewrite} and representations \eqref{EQ:repre} and \eqref{EQ:repre-d}, we write
\begin{equation}\label{EQ:start}
\begin{split}
\|\nabla u(t) \|_{L^p(\Omega)}
& \ge \left\| -(r+1)^{-2+\frac{2}{p}} v(t) + (r+1)^{-1+\frac{2}{p}} \partial_r v(t) \right\|_{L^p(0, \infty)}\\
& =  (4 \pi t)^{-\frac{1}{2}} \left\| (r+1)^{-1+\frac{2}{p}} \int_0^\infty K(t,r,s) g(s)\,ds\right\|_{L^p(0, \infty)},
\end{split}
\end{equation}
where
\begin{equation*}
K(t,r,s)
= \left[ \left\{-(r+1)^{-1}  - \frac{r-s}{2t}\right\} e^{-\frac{(r-s)^2}{4t}} + \left\{(r+1)^{-1} +\frac{r+s}{2t} \right\} e^{-\frac{(r+s)^2}{4t}} \right].
\end{equation*}
Again we take $t=t_m=m^2$ and denote by $u_{m}$ and $v_{m}$ the solutions to \eqref{EQ:DP} and \eqref{EQ:DP-rad} with initial data $f_{m}$ in \eqref{EQ:f_m} and $g_{m}$ in \eqref{EQ:g_m}, respectively.

To begin with, we prove the following:
For sufficiently large $m\in\mathbb N$, there exists a constant $C>0$, independent of $m$, such that
\begin{equation}\label{EQ:K-below}
K(m^2,r,s) \ge \frac{C}{m}
\end{equation}
for any $10\le r \le m^{1/4}$ and $m\le s\le 2m$.
Writing
\[
e^{-\frac{(r-s)^2}{4t}} = e^{-\frac{s^2}{4m^2} + \frac{r}{2m} - \frac{r^2}{4m^2}}
= e^{-\frac{s^2}{4m^2}} \left\{ 1 + \frac{r}{2m} + \mathrm{O} \left(\frac{r^2}{m^2} \right) \right\},
\]
\[
e^{-\frac{(r+s)^2}{4t}} = e^{-\frac{s^2}{4m^2} - \frac{r}{2m} - \frac{r^2}{4m^2}}
= e^{-\frac{s^2}{4m^2}} \left\{ 1 - \frac{r}{2m} + \mathrm{O} \left(\frac{r^2}{m^2} \right) \right\},
\]
we calculate
\begin{equation}\label{EQ:K}
\begin{split}
K(m^2,r,s)
& = e^{-\frac{s^2}{4m^2}} \left[ \left\{-(r+1)^{-1}  - \frac{r-s}{2m^2}\right\} \left\{ 1 + \frac{r}{2m} + \mathrm{O} \left(\frac{r^2}{m^2} \right) \right\}\right. \\
&  \qquad\qquad\qquad \left. + \left\{(r+1)^{-1} +\frac{r+s}{2m^2} \right\} \left\{ 1 - \frac{r}{2m} + \mathrm{O} \left(\frac{r^2}{m^2} \right) \right\} \right]\\
& =  e^{-\frac{s^2}{4m^2}} \left\{ \frac{s}{m^2} -(r+1)^{-1}\frac{r}{m}- \frac{r^2}{2m^3} +\mathrm{O} \left(\frac{r^2}{m^2} \right) \right\}\\
& \ge e^{-1} \left\{ \frac{1}{11m}
+\mathrm{O} \left(\frac{r^2}{m^2} \right) \right\},
\end{split}
\end{equation}
where we used in the last step
\[
\frac{s}{m^2} -(r+1)^{-1}\frac{r}{m} \ge \frac{1}{m} - \frac{10}{11m} \ge \frac{1}{11m}
\]
for $10\le r \le m^{1/4}$ and $m\le s\le 2m$.
Since we can neglect the remainder terms in \eqref{EQ:K} if $m$ is sufficiently large,
we obtain \eqref{EQ:K-below}.

Let us turn to estimate form below of $L^p$-norm of $\nabla u_m(t_m)$.
By combining \eqref{EQ:start} and \eqref{EQ:K-below}, we estimate
\begin{equation} \label{EQ:u_m4}
\begin{split}
\|\nabla u_m(t_m) \|_{L^p(\Omega)}
& \ge C m^{-2} \left\| (r+1)^{-1+\frac{2}{p}} \int_m^{2m}  C_m\,ds
\right\|_{L^p(10,m^{\frac{1}{4}})}\\
& =
C\cdot C_m m^{-1} \left\| (r+1)^{-1+\frac{2}{p}}\right\|_{L^p(10,m^{\frac{1}{4}})}\\
& \ge C\cdot C_m m^{-1}
\end{split}
\end{equation}
for sufficiently large $m\in\mathbb N$.
Noting from \eqref{EQ:C_m} that
\[
C_m m^{-1} \sim m^{1-\frac{3}{p}} m^{-1}=m^{-\frac{3}{p}}
\]
as $m\to \infty$, we conclude from \eqref{EQ:u_m4} that
\[
\|\nabla u_m(t_m) \|_{L^p(\Omega)} \ge C m^{-\frac{3}{p}} = C t_m^{-\frac{3}{2p}}
\]
for sufficiently large $m\in\mathbb N$, where the constant $C>0$ is independent of $m$.
This proves that
\[
\limsup_{m\to\infty}\frac{t_m^{\frac{3}{2p}}\|\nabla u_{m}(t_m)\|_{L^p(\Omega)}}{\|f_m\|_{L^p(\Omega)}} >0,
\]
since $\|f_m\|_{L^p(\Omega)}=1$ by \eqref{EQ:f_m}.
Thus the optimality for $3< p\le \infty$ is proved. The proof of Theorem \ref{thm:optimal} is finished.
\end{proof}

\appendix
\section{}
\label{App:A}
In this appendix, we show the estimate \eqref{EQ:opti-infi}:
\[
\|\nabla e^{-tH}\|_{L^\infty(\Omega)\to L^\infty(\Omega)} \ge C
\]
for any $t>1$, when $n\ge 3$ and $\Omega$ is the exterior domain of a compact connected set with $C^{1,1}$ boundary.
The estimate \eqref{EQ:opti-infi} follows from the known result:
There exists a constant $C>0$ such that
\begin{equation}\label{EQ:opti-1-infi}
\|\nabla e^{-tH}\|_{L^1(\Omega)\to L^\infty(\Omega)} \ge C t^{-\frac{n}{2}}
\end{equation}
for any $t>1$ (see Section 1 in \cite{IshKab-2007} and also \cite{Zhang-2003}).
In fact, we suppose that
\begin{equation}\label{EQ:supp}
\|\nabla  e^{-tH} \|_{L^\infty(\Omega)\to L^\infty(\Omega)} \le C(t)
\end{equation}
for any $t>1$, where $C(t)\to \infty$ as $t\to \infty$.
Then we deduce from \eqref{EQ:supp} and Lemma \ref{lem:LpLq} that
\[
\|\nabla  e^{-tH} f\|_{L^\infty(\Omega)} \le C(t) \|e^{-\frac{t}{2}H}f\|_{L^\infty(\Omega)}
\le C\cdot C(t) t^{-\frac{n}{2}}\|f\|_{L^1(\Omega)}
\]
for any $t>1$ and $f\in L^1(\Omega)$.
This contradicts the fact \eqref{EQ:opti-1-infi}. Thus \eqref{EQ:opti-infi} is true.

\section{}
\label{App:B}
In this appendix we prepare two fundamental inequalities.
The first one is the special case of the Gagliardo - Nirenberg inequality (see \cite{Gag-1959}, \cite{Nir-1959}).

\begin{lem}
\label{lem:GN}
Let $\Omega$ be a bounded domain in $\mathbb R^n$ having the cone property.
Then, for any $1 < p < \infty$ and multi-index $\alpha$ with $|\alpha|=2$, there exist constants $C_1, C_2 >0$ such that
\[
\|\nabla f\|_{L^p(\Omega)} \le C_1
\|D^\alpha f\|_{L^p(\Omega)}^{\frac{1}{2}} \|f\|_{L^p(\Omega)}^{\frac{1}{2}}
+ C_2 \|f\|_{L^p(\Omega)}
\]
for any $f\in W^{2,p}(\Omega)$.
\end{lem}

The second one is the global $W^{2,p}$-estimate (see Theorem 9.14 in \cite{GilTru_2001}).

\begin{lem}\label{lem:B}
Let $\Omega$ be a domain in $\mathbb R^n$ with $C^{1,1}$ boundary.
Then, for any $1 < p < \infty$, there exists a constant $C>0$ such that
\[
\|f\|_{W^{2,p}(\Omega)} \le C \left(\|\Delta f \|_{L^p(\Omega)} + \| f \|_{L^p(\Omega)}\right)
\]
for any $f \in W^{2,p}(\Omega) \cap W^{1,p}_0(\Omega)$.
\end{lem}


\end{document}